\documentclass[amssymb]{article}
\def\XN{{ X_0(N) }}
\def\JN{{ J_0(N) }}

\def\Hom{{ \mbox{Hom} }}
\def\Aut{{ \mbox{Aut} }}

\def\T{{ {\bf T} }}
\def\Tr{{ \mbox{Tr} }}
\def\1ox{{ \Omega^1_{\scriptstyle{X}} }}
\def\2ox{{ \Omega^2_{\scriptstyle{X}} }}

\def\ok1{{ \Omega^1_K }}
\def\ok2{{ \Omega^2_K }}

\def\P{{ \mathbb{P} }}

\def\ra{{ \rightarrow }}

\def\hra{{ \hookrightarrow }}

\def\F{{ {\bf F} }}

\def\Q{{ {\bf Q} }}

\def\C{{ {\bf C} }}

\def\8{{ {\infty } }}

\def\G{{ \Gamma }}
\def\Gal{{ \mbox{Gal} }}

\def\^{{ ^{\wedge} }}

\def\hZ{{ \hat{ \bf Z} }}
\def\Z{{ {\bf Z } }}

\def\tor{{ \mbox{tor} }}

\newtheorem{thm}{Theorem}

\newtheorem{Def}{Definition}
\newtheorem{lem}{Lemma}

\usepackage{amsmath}
\usepackage{amsfonts}
\usepackage{amscd}
\usepackage{amssymb}

\title{ Torsion
points on modular curves and Galois Theory}
\author{ Kenneth A. Ribet and Minhyong Kim}

\begin{document}

\maketitle

In elementary terms,
the arithmetic theory of a curve $X$ is concerned with
solutions to a
 geometrically irreducible polynomial equation in two variables:
$$ f(x,y)=0\ \ \ \ \ (*)$$
In contrast to the geometric theory, where the
different kinds of number pairs $(x,y)$ that can occur
as solutions are viewed as homogeneous, the arithmetic study classifies more carefully the structure of
solutions of specific  type. That is, one tries to
understand  the solutions to the equation (*)
where $(x,y)$ are constrained to lie in
some arithmetically defined set.
 One common case is that of
rational solutions or, more generally, the
case of solutions where $x,y$ are
constrained to lie inside  a fixed number field $F$.
For example, when $f(x,y)$ has genus 1  (that is, the smooth
points of the complex solution set form a genus one Riemann surface
with punctures), the Mordell-Weil theorem says the solution
set, in conjunction  with a few additional points, acquires the
natural structure of a finitely generated abelian group.
For another example, when the genus is greater than one, Faltings \cite{Fa}
 proved
that the solution set is finite. In both cases, one derives
finite-type structures for  the solution set from
finiteness contraints on the `type' of the solution.
A theorem of Ihara-Serre-Tate (\cite{La1}, theorem 8.6.1)
 gives an example of finiteness
theorems deriving from a different kind of arithmetic
constraint. Namely, one considers solutions that are
roots of unity of arbitrary order. Then as soon as the
genus is at least one, there are again only finitely
many solutions to (*). It is interesting to note that
in this case, the constraint in question is
of `group-type' in contrast to the `field-type'
constraint of the other two examples.

A conjecture of Manin and Mumford as proved by
Raynaud \cite{Ra1} deals with the
projective case of this theorem.
What is meant by the projective case? In the Ihara-Serre-Tate
theorem, one can view the curve $X$ as essentially lying in
the affine torus $\C^*\times \C^*$ and the assertion
is that $X$ has a finite intersection with the
torsion points of the torus. Now, when $X$ is a projective
smooth curve of genus at least two,
it has an essentially canonical embedding
into a group variety $J=J(X)$, the Jacobian of $X$. Raynaud's
theorem 
 states that the intersection between $X$ and $J_{tor}$,
the torsion subgroup of $J$, is finite. It should also be
noted that Raynaud generalizes this
statement considerably to include refined statements
about intersections between  subvarieties
of abelian varieties and division points of
finitely generated groups \cite{Ra2}, while
a common generalization of the
projective and affine case concerned with subvarieties of
semi-abelian varieties has been found by Hindry \cite{Hi}.

On the other hand, Coleman \cite{Co1} \cite{Co2} \cite{Co3}
has stressed
the importance of  being able to determine explicitly
the finite set occuring in Raynaud's theorem
for specific curves.
This program has been carried out with some success, most notably 
in the case of Fermat curves, due to
Coleman, Tamagawa, and Tzermias \cite{CTT},
and the modular curves $X_0(N)$ for
$N$ prime, due to Baker \cite{Ba} and Tamagawa.
A new proof for the modular curve case
was given by Ken Ribet using a refined
analysis of the Eisenstein torsion in $J_0(N)$
and this paper is devoted to an exposition of
this proof. It is similar in many ways to the second proof of
\cite{Ba}  except for conceptual simplifications arising
from  systematic use of the notion of an `almost rational
torsion point.' In particular, a complete computation of
these points is given for $\JN$, and Lemma 1 makes clear
how the main theorem hinges upon this notion.
The
result in question was first conjectured
by
Coleman, Kaskel, and Ribet \cite{CKR} and we go on to
describe the statement. As mentioned,
we will always be interested in prime $N$
that are $ \geq 23$ (which occurs iff $X_0(N)$ has genus
$\geq 2$).  $X_0(N)$ has two cusps
corresponding to the orbits of
$0$ and $\infty$ in the extended upper half-plane,
and we will use the latter, again denoted
$\infty$, to embed $X_0(N)$ into its
Jacobian $i:X_0(N) \hra J_0(N)$. That is, a point
$P \in X_0(N)$  maps to the class of the
divisor $[P]-[\infty]$. In the following,
we will suppress the embedding $i$ from
the notation or leave it in according to
convenience. By a
 theorem of Manin and Drinfeld \cite{Dr}, the  cusp 
$0$ is a torsion point under this
embedding. Another way for a torsion point to
arise is as follows: the curve
$X_0(N)$ has an involution $w$ which switches
$0$ and $\infty$, that is, 
$0=w(\infty)$. Denote by $\XN^+$ the curve
obtained as the quotient of $X_0(N)$ by the
action of this involution. Now, it can happen
that $\XN^+$ is a curve of genus zero,
in which case $X_0(N)$ is a hyperelliptic curve.
Let $$f: \XN \ra \XN^+$$ be the quotient map
and let $P \in \XN$ be a Weierstrass  point.
The inverse image divisors
of any two points are rationally equivalent, since
$\XN^+\simeq \P^1$. In particular,
$2[P]\sim [\infty]+[0]$. Thus, $2i(P)=i(0)$
and $i(P)$ is a torsion point.  
According to Ogg \cite{Og}, the values of $N$ for which
$\XN$ is hyperelliptic are
$23,29,31,37,41,47,59,71$. In the case  $N=37$
the hyperelliptic involution $h$ is
different from $w$. That is $X_0(37)$ is
hyperelliptic even though  $X_0(37)^+$ is
not of genus zero. It was shown by
Mazur and Swinnerton-Dyer \cite{MS} that 
 $[\infty]-[h(\infty)]$ is of infinite order
in $J_0(37)$. From this  it is an easy exercise to deduce  that
the Weierstrass points are not torsion in this case.
That is, Weierstrass torsion points occur
only when $\XN^+$ is of genus zero.
Thus, we will have completely
determined the torsion points as soon as we have
found the non-Weierstrass ones.

 The conjecture of
Coleman, Kaskel, and Ribet as proved by
Baker and Tamagawa says, in fact,  the following:
\begin{thm} (Baker, Tamagawa)
$$[X_0(N)(\overline{\Q})-\mbox{(Weierstrass points)}] \cap J_{\tor}=\{0, \infty\}.$$
\end{thm}

\section{Almost rational torsion points}
Lang's original suggestion \cite{La2} was to prove
the Manin-Mumford conjecture itself using
Galois theory. 
Let's suppose given a curve $C$ embedded
in an abelian variety $A$ over the complex numbers. The data is
defined over some field $K$ finitely generated
over the rationals, and hence, the torsion
points of $A$ will admit an action of
the Galois group $G=\Gal(\overline{K}/K)$. This
action induces a representation
$$\rho:G \ra \Aut( \hat{T}A)$$
where $\hat{T}A$ denotes the adelic Tate module
of $A$. Lang's conjecture concerns the
intersection between the image $\rho(G)$
 of $G$  and the group of homotheties $\hat{\Z}^* \subset \Aut (\hat{T}A)$.
He conjectured that $\rho(G) \cap \hat{\Z}^*$ is
of finite index in $\hat{\Z}^*$. The Manin-Mumford
conjecture follows from this by an elementary
argument.

Although Lang's conjecture is still unproven,
Serre proved a weaker version in his College
de France lectures 85-86 \cite{Se}. That is, he proved
that
$\hZ^*/\rho(G)\cap \hZ^*$ is of finite
exponent.
Using Serre's result Ribet manages to give a very elegant
proof of the Manin-Mumford conjecture.

In this proof  crucial use is made of
the notion of an `almost rational' torsion points, which we
will abbreviate to a.r.t.:

\begin{Def}
Let $A$ be an abelian variety over a field $k$.
A point $p\in A(\bar{k})$ is called almost
rational (a.r.) if 
$$\sigma (p)-p=p-\tau (p) \Rightarrow p=\sigma (p)=\tau(p)$$
for all $\sigma , \tau \in \Gal (\bar{k}/k)$.
\end{Def}
Here are a few elementary facts that follow directly
from the
definition:

-Rational points are almost rational.

-A Galois conjugate of an a.r. point is a.r.

-If $P$ is almost rational and $2\sigma (P)-2P=0$
then $\sigma (P)=P$.

Even after verifying these facts,
the definition is not likely to be very intuitive,
so it is probably best to see right away a concrete result that
uses it.

\begin{lem}
Let $X$ be a curve of genus at least 2 embedded in
its Jacobian $J$ via a rational point $p_0$.
Then $$X=X_{a.r.}\cup (\mbox {Weierstrass points})$$
\end{lem}

Thus, we get an inclusion $X_{tor}-(\mbox{w.p.})\subset J_{a.r.t.}$
reducing the Baker-Tamagawa theorem to  a study of
$J_{a.r.t.}$ for $J_0(N)$.

\smallskip

{\em Proof of Lemma.}
Suppose $[P]-[P_0]$ is not almost rational. Then there are
$\sigma$ and $\tau$ in the Galois group such that
$[\sigma (P)]-[P ]\sim  [P]-[\tau (P)]$ as divisors and neither
are equivalent to zero. Thus,
$2[P]- [\sigma (P)]-[\tau (P)] \sim 0$, meaning we can find
a rational function with a pole of order two at $P$.
That is, $P$ is a Weierstrass point.

We will investigate this notion extensively in the specific
context of modular curves in order to prove the Baker-Tamagawa
theorem. In the meanwhile, we  outline how
to  deduce the Manin-Mumford conjecture from
Serre's result.
In fact, Manin-Mumford obviously follows  from
Lemma 1 and the following theorem, whose proof will
occupy us to the end of this section.
\begin{thm}
Let $A$ be an abelian variety over a finitely generated
field $k$. Then $A_{a.r.t}$ is
finite.
\end{thm}

In the course of the proof, we will need the following simple
\begin{lem}
For each $e\geq 1$, we can find $C(e)>0$ such that
for any $m> C(e)$, there exist $x,y \in ((\Z/m\Z)^*)^e$ with
$x\neq 1, y\neq 1$ and $x+y=2$.
\end{lem}

{\em Proof.}
First note that if $m=\prod p^{n_p}$, then by the Chinese
remainder theorem,
one need prove the existence of $x,y$ for just one of
the $\Z/p^{n_p}\Z$ and set the modulus for the other
factors to be 1. Also, by setting $C(e)$ sufficiently
large, we can make sure that there is at least
one prime power factor $p^{n}\geq A(e)$, where $A(e)$
is the maximum of $e^4$ and 1+the biggest prime $l$ such that
$x^e+y^e=2$ has at most $e^2+2e$ solutions in $\F_l$.
Such an $l$ clearly exists by elementary counting when
$e$ is  1 or 2 and by the Weil bounds when $e\geq 3$.

In the case $n\geq 2$  write $e=up^k$ where $u$ is
relatively prime to $p$. Now put $x=1+ep^{n-k-1}$
and $y=1-ep^{n-k-1}$ and note that $p^n\geq e^5=u^3p^{5k}$
implies that $k\leq \lfloor n/5$ and $k=0$ for
$n\leq 4$, so that, in any case, $x$ and $y$ are
both units in $\Z/p^n$. Clearly $x,y \neq 1$ (mod $p^n$)
but $x+y=2$ (mod $p^n$). It is also easily checked that $x=(1+p^{n-k-1})^e$
and $y=(1-p^{n-k-1})^e$ (mod $p^n$).
Next suppose $n=1$. Then we are looking for solutions
to $x^e+y^e=2$ in $\F_p$ such that neither  $x^e$ nor $y^e$
are $0$ or 1. We are done by counting the number of points mod $p$.

It is easy to sharpen the proof slightly and take $C(e)=3$ if
$e=1$.
\smallskip

{\em Proof of theorem.} According to Serre,
if we consider the action $\rho :G \ra \Aut (\hat{T}A)$
of the Galois group
on the adelic Tate module, 
$\hZ^*/\rho(G)\cap \hZ^*$  has finite exponent $e$.
We claim that if $P$ is a torsion
point of order $ m> C(e)$, the $P$ is not a.r.
To see this, let $x,y \in ((\Z/m\Z)^*)^e$ satisfy
the conditions of the proposition. Find $\sigma, \tau \in G$
such that $\sigma \mapsto x$ and $\tau \mapsto y$ as operators on $A[m]$.
Then we have
$$\sigma (P)+ \tau (P)=2P\Rightarrow \sigma (P)-P=P-\tau (P)$$
but $\sigma (P)-P=(x-1)P\neq 0$. That is, $P$
is not almost rational.

\smallskip

Thus, the finiteness of a.r.t.~ points
follows from very general considerations.
To prove the target theorem in the case of modular
curves, we will end up needing a very precise understanding
of the a.r.t.~ points for modular Jacobians, in particular,
their relation to other canonically defined subgroups
with special  Galois-theoretic properties. We will
review the relevant facts in the next section.

We close this section with a few lemmas for use in the proof
of the main theorem.

\begin{lem}
Let $A/\Q$ be an abelian variety and suppose $P\in A[n]$, $n>3$,
is a cyclotomic point, i.e., $\sigma (P)=\chi_n (\sigma )P$
for all $\sigma \in \Gal (\overline{\Q}/\Q)$,
where $\chi_n$ is the mod $n$ cyclotomic character.
Then $P$ is not a.r.
\end{lem}
\smallskip

{\em Proof.}
As noted above, it is easy to see that if $n>3$, then there exist
$s,t \in (\Z/n\Z)^*$ such that $s\neq 1, t\neq 1 $
and $s+t=2$. Find $\sigma , \tau$ such that
$\chi_n(\sigma )=s$ and $\chi_n(\tau)=t$. Then
$\sigma (P)+\tau (P)=2P$ but $\sigma (P)-P=sP-P \neq 0$.
So $P$ is not a.r.
\smallskip

\begin{lem}
Let $A$ be deinfed over a number field  $k$.
 Let $v$ be prime of $k$ and
assume $A$ has semi-stable reduction at $v$. Let $P\in A_{a.r.t}$
have order prime to $v$. Then $k(P)$ is unramified at $v$.
\end{lem}
\smallskip

{\em Proof.} Let $\sigma \in I_v$, an inertia group at
$v$. According to Grothendieck (\cite{Gr}, see also following section), 
the action of $I_v$
on prime to $v$ torsion is two-step unipotent. So
$$\begin{array}{ccc}
(\sigma -1)^2 P=0 & \Rightarrow & \sigma ^2 P-2\sigma P+P=0\\
 & \Rightarrow & \sigma P+\sigma^{-1}P= 2P \\
&\Rightarrow & \sigma (P)=P
\end{array}$$
the last implication following from the assumption
that $P$ is a.r. Therefore, $I_v$ acts trivially on $P$.

\section{Background on Modular curves}
In this section, we summarize
 the facts we need from the  theory of modular  curves,
especially
results about the Galois representations associated to their
Jacobians. (See \cite{Ma} and references therein for a general
overview.)

Recall that the modular curve $X_0(N)$ is the projective
smooth model of the modular curve $Y_0(N)$ which
parametrizes pairs $(E,C)$, where $E$ is an elliptic
curve and $C$ is a cyclic subgroup of order $N$.
$Y_0(N)$ and $X_0(N)$ are defined over
$\Q$, and
over the complex numbers, we have
$$Y_0(N)(\C)=H /\G_0(N)$$
while $$X_0(N)(\C)=[H \cup \P^1(\Q)]/\G_0(N)$$
When $N$ is prime, which is the case that will
concern us, $\G_0(N)$ has two orbits on
$\P^1(\Q)$, the orbits of 0 and $\infty$.
We will denote by the same symbols the
corresponding points on $X_0(N)$.
We denote by $\JN$ the Jacobian of $\XN$,
which  parametrizes  divisor
classes of degree zero on $\XN$.
The Abel-Jacobi embedding $\XN \hra \JN$ with respect to
the point $\infty$ is described at the
level of points by sending a point $P$ to
the class of the divisor 
$[P]-[\infty]$. We will use this to
 identify $\XN$ with its image and think
of it as a subvariety of $\JN$.
The Manin-Drinfeld
theorem says that
$[0]-[\infty]$ generates a finite subgroup
$C$ of $\JN$ which we call the {\em cuspidal
subgroup}. We will denote by $n$ the order of
$C$, which is equal to the numerator of
$(N-1)/12$ (\cite{Ma} p. 99). 

Another important subgroup is the
Shimura subgroup $\Sigma$ of $J_0(N)$ defined as follows.
There is a map $X_1(N) \ra X_0(N)$ of degree $(N-1)/2$
 from the compactification
$X_1(N)$ of the modular curve $Y_1(N)$ which
parametrizes pairs $(E,P)$, where $E$ is an elliptic curve
and $P$ is a point of order $N$. On the points of
$Y_1(N)$ this map simply takes
$(E,P)$ to $(E,<P>)$, $<P>$ being the subgroup
generated by $P$.  This gives rise to a map $X_2(N)\ra \XN$
which is the maximal \'{e}tale intermediate covering
to $X_1(N) \ra \XN$. 
Thus we get a map $J_2(N) \ra J_0(N)$, where
$J_2(N)$ is the Jacobian of $X_2(N)$.
$\Sigma $ is simply the kernel of the dual map. Thus,
the points of $\Sigma$ correspond to line bundles
of degree zero on $X_0(N)$ which become trivial
when pulled back to $X_2(N)$. It has order $n$ and is
 isomorphic to $\mu_n$ as a Galois
module (\cite{Ma} p.99).

The modular Jacobians admit an action of the
algebra $\T$ of Hecke operators (\cite{Ma}, section II.6). This
is the $\Z-$algebra of endomorphisms generated
by the correspondences  $T_l$ for
each prime $l\neq N$ and the Atkin-Lehner
involution $w_N$. They are defined on
points of $Y_0(N)$ by the formula
$$T_l: (E,C)\mapsto \Sigma_{C'} (E/C', (C+C')/C')$$
where $C'$ runs over the cyclic subgroups of $E$
of order $l$ and
$$w_N: (E,C)\mapsto (E/C, E[N]/C).$$

The Eisenstein ideal $I$
of $\T$ is the ideal generated by $T_l-(l+1)$
for $l\neq N$ and $1+w_N$ (\cite{Ma} p.95). Of particular importance
will be the structure of
the subgroup $\JN [I]\subset \JN$ annihilated by $I$. 
\smallskip

We now list the main difficult facts we will use:

\medskip
(0) $\T/I \simeq \Z/n$ (\cite{Ma}, Prop. II.9.7). So if a maximal ideal
$m$ is `Eisenstein', i.e., contains $I$, then
$\T/m$ has characteristic $l$ dividing $n$.
\medskip

(1)
$\JN [I]=C\oplus \Sigma$ if $n$ is odd while
$\JN [I]$ contains $C +\Sigma$ as a subgroup of
index two and $C\cap \Sigma =C[2]=\Sigma [2]$
if $n$ is even. This follows from the fact that $C+\Sigma$
is contained in $\JN[I]$ and that $\JN [I]$ is
free of rank two over $\T /I$. (See \cite{Ma}, sections II.16-18, and 
Prop. II.11.11
together with the explanation in \cite{Ri2}, section 3.)

\medskip
(2) We will  need some detailed facts about the action
of the Galois group $G=\Gal (\bar{\Q}/\Q)$ on  the torsion points
of $\JN$. One analyzes these representations
by breaking them up into simple $\T[G]$-modules.
Such simple modules are associated to maximal
ideals $m$ inside the Hecke algebra $\T$.
In fact, for each $m$ there is a  two-dimensional semi-simple
representation, unique up to isomorphism,
$$\rho_m:G \ra GL_2(\T/m)$$
characterized by the properties (\cite{Ri1}, section 5):

-$\rho_m$ is unramified outside $N$ and $l$, where
$l=m\cap \Z$.

For $p\neq N,l$, the Frobenii $Fr_p$ satisfy

-$\Tr (\rho_m(Fr_p))=T_p$ (mod $m$) 

-and $\det (\rho_m(Fr_p))=p.$ 

Furthermore, one knows that $\rho_m$ is irreducible
if $m$ is non-Eisenstein, i.e., when $m$ does not contain
the Eisenstein ideal $I$, and if $I\subset m$, then
$\rho_m$
is isomorphic to $\Z/l\oplus \mu_l$ (\cite{Ma} Prop. 14.1 and 14.2).  

\medskip
(3) Concerning the action of $I_N$, the inertia
group at $N$, on the torsion of $\JN$, one has
Grothendieck's exact sequence (\cite{Gr} 11.6 and 11.7)
$$0\ra \Hom (X, \mu_r) \ra \JN [r] \ra X/rX \ra 0$$
for any $r$, where $X$ is the character group of the
toric part of the reduction of $\JN$ mod $N$.
This implies, for example, that the action is
2-step unipotent if $r$ is prime to $N$.
One notes also that even when $N |r$, the
first and last terms are finite, in that they
extend to finite flat group schemes over $\Z_N$.

\medskip
(4) On the other hand, a theorem of Ribet (\cite{Ri2} Prop. 2.2) 
addresses fine behaviour of
 $\rho_m$ at $N$ for $m$ non-Eisenstein. It says
that $\rho_m$ is not finite at $N$ if $m|N$ and that it
is ramified at $N$ if $m$ is prime to $N$.
This is an instance of the `level-lowering'
theorem \cite{Ri1}, together with a result of Tate on
mod 2 representations unramified outside 2 \cite{Ta}. 

For $m|N$, $\rho_m$ occurs in 
$\JN[N]$, so as an $I_N$ module, it fits into
an exact sequence
$$0\ra \mu_N \ra \rho_m \ra \Z/N \ra 0$$
which is non-split, since the existence of
a splitting would
imply finiteness for $\rho_m$. So we draw the 
conclusion that $\rho_m (I_N)$ is
non-abelian in this case.

\medskip
(5) 
A theorem of Ribet says
that $\JN [I]$ is exactly the set of
 torsion in $\JN $ that is unramified at $N$.
(\cite{Ri2} Prop. 3.1, 3.2) That $\JN[I]$ is
unramified at $N$ for $n$ odd follows obviously from
$\JN[I]=C\oplus \Sigma$. When $n$ is even, one still
gets an isomorphism $$\JN[I]\simeq\Hom (X/IX, \mu_n)\oplus \Sigma.$$
That is,
 the two groups on right hand side inject into the left
by (1) and (3) and the images do not intersect \cite{Ma}, (Prop. II.11.9).
But they also have the same order by (1) and the argument
of \cite{Ri2} theorem 2.3 showing that $X/IX$ is cyclic.

To go the other way, given an unramified torsion point
$P\in \JN$, one uses (4) to conclude that the simple
constituents of the module $M:=\T[G]P+\JN[I]$ all come from Eisenstein
primes, and therefore, are of the form $\Z/l$ or $\mu_l$ for
$l|n$. So the constituents are all annihilated by $I$. 
 It is easy to
see then that $M$ itself is of the form
$$0\ra S \ra M \ra Q \ra 0$$
where $Q$ is constant and $S$ is of $\mu$-type. But $\Sigma$
is the maximal $\mu$-type group in $\JN$ (\cite{Ma}, theorem 2)
so $S= \Sigma$. Now, reduction mod $N$ and the isomorphism
between $\Sigma$ and the component group of $\JN$ mod $N$ gives
us a splitting of this exact sequence. So one need only show
that $I$ annihilates $Q$. The Eichler-Shimura relation say that
$$T_l\cong Fr_l+lFr_l^t$$ (mod $l$, $l\neq N$), and therefore, the
constant group 
$Q$ is annihilated by $\eta_l=T_l-(1+l)$ for $l\neq N, (l,n)=1$.
(The order of $Q$ divides some power of $n$, so reduction mod $l$
is injective on $Q$ for $l$ prime to $n$.)
 
To show that it is
also annihilated by all of $I$, we decompose into $m$-primary
factors for Eisenstein primes $m$
(which is possible since $Q$ is annihilated by some power of $I$)
and then show that each factor is annihilated using
local principality of the Eisenstein ideal (\cite{Ma}, theorem II.18.10).

\medskip
(6) According to a theorem
of Ribet (\cite{Ri2} Theorem I.7),  the field $\Q(\JN [I])$ generated by
the Eisenstein torsion is $\Q(\mu_{2n})$ while
$\Q(C,\Sigma)=\Q (\mu_n)$. The proof of the first fact follows
from a careful study of $\JN[I]$, but appears a bit too elaborate
to summarize in a few words. On the other hand, note that
for $n$ odd, the first fact follows easily from the second.

\medskip

(7) Finally, it is explained by Coleman-Kaskel-Ribet \cite{CKR}
that Mazur's theorems imply the useful fact that
$\XN \cap C =\{0,\infty \}$. 
For $N\neq 37, 43, 67, 163$, it is
an obvious consequence of the fact that the cusps are
the only rational points of $\XN$.  The remaining cases
can be treated by more elementary arguments.

\section{The theorem of Baker-Tamagawa}

The main result which provides the key is the following
\begin{thm}
$\JN_{a.r.t}=C\oplus \Sigma [3]$
\end{thm}
This detailed knowledge is what makes it possible to
determine the torsion points on $\XN$ so explicitly.

Let us first show how theorem 3 implies the theorem
of Baker and Tamagawa.

This implication divides into two cases. Recall the
curve $\XN^+$ obtained as the quotient of
$\XN$ by the action of $w$, the Atkin-Lehner involution.
The first case is when $\XN^+$ has positive genus.
Then the projection $$f:\XN \ra \XN^+$$
induces   a commutative diagram:
$$\begin{array}{ccc}
\XN & \hra & \JN \\
\downarrow & & \downarrow \\
\XN^+ & \hra & \JN^+
\end{array}
$$
where $\JN^+$ denotes the Jacobian of $\XN^+$. 
According to the
theorem,
 $$\JN_{a.r.t.}\subset \JN[I] \subset \JN[1+w]$$
Now, if $D$ is a degree zero divisor on $\XN$,
then $$D+wD=f^*f_*(D).$$ So if $D+wD\sim 0$, then
the class of $f_*(D)$ is in the kernel of
$$f^*:\JN^+\ra \JN .$$ But since $w$ has a fixed
point, this map is injective. Thus,
$$\JN[1+w] \ra 0 \in \JN^+,$$ and therefore,
$$\JN_{a.r.t}\ra 0.$$ But this implies that
$$\XN_{tor} \ra \infty  \in \XN^+$$
and hence that $\XN_{tor}=\{0,\infty \}$
as desired.

The second case is when $\XN^+$ fails to have positive genus,
that is,
when $N=23,29,31,41,47,59,71$. In this
case, $N$ is not congruent to 1 mod 9 which in turn
implies that  3 does not divide $n$. Therefore,
by  theorem 2 $\JN_{a.r.t.}=C$, and  we get
$$\XN_{tor}\subset \XN \cap C=\{0,\infty \}$$ again.

\smallskip
So it remains to prove the structure theorem for
$\JN_{a.r.t.}$.

\smallskip
We wish to show first that $\JN_{a.r.t.}\subset \JN [I]$,
which is the hard part of the proof. This is achieved by
proving that the points in  $\JN_{a.r.t.}$
are unramified at $N$, and using Ribet's theorem
identifying such points with $\JN [I]$.

 To prove that $\JN_{a.r.t}$
consists of points unramified over $N$ it suffices
to show that the points have order prime to $N$ (Lemma 4).
So let $P\in \JN_{a.r.t.}$ and analyze the module
$M:=\T[G] P$ by breaking it into its simple constituents,
the possibilities for which we described in the previous
section.
Let $r$ be the order of $P$. Thus, we have
$M\subset \JN[r]$. 

In order to see that $\JN_{a.r.t}\subset \JN [I]$,
recall from the previous section that as an $I_N$ module, $\JN[r]$
fits into an exact sequence
$$0\ra \Hom (X,\mu_r) \ra \JN[r]\ra X/rX \ra 0$$
Therefore $$I'_N:=\mbox{Ker}( \chi_r :I_N \ra  (\Z/r)^*)$$
acts on $\JN[r]$ by two-step unipotent transformations.
But this implies by the argument of Lemma 4 that
$\sigma (P)=P$ for all $\sigma \in I'_N$.
The same argument also applies to the conjugates of $P$
since they are also a.r.t.
Therefore, $I'_N$
acts trivially on  $M$. That is, $I_N$ acts through the
quotient $I_N/I'_N \hra (\Z/r)^*$ on $M$ and
all its constituents. From this, we see 
that $\rho_m$ for $m|N$
is ruled out as a simple factor (since $\rho_m(I_N)$
is
non-abelian in that case) leaving $\Z/l$, $\mu_l$,
and $\rho_m$, for $m$ not dividing $N$, as possibilities.
Since the one-dimensional factors only occur in the Eisenstein
case, we get $l|n$ and therefore, $l$ is relatively prime to $N$.
We conclude that $M$ must have order prime to $N$, and hence, so
must $P$. Therefore, $P \in \JN[I]$ as desired.

In fact, we claim that $P\in \Sigma +C$.
For if $P\notin \Sigma +C$ (which occurs only
when $n$ is even), $P$ must generate $\JN[I]/(\Sigma +C)$,
so by fact (6) of the previous section, 
 we must have $\Q (P,\Sigma, C)=\Q (\mu_{2n})$.
 Also, $\Q (\Sigma, C)= \Q (\mu_n)$.
Therefore,
 we can find $\sigma \in G$ such that $\sigma (P)-P \neq 0$
and $\sigma $ acts trivially on $C+\Sigma$. But we have
$2P \in \Sigma +C$, so that $\sigma (2P)-2P=0$.
This contradicts the assumption that $P$ is a.r.
by our remark following the definition of a.r.

So we have $P \in \Sigma +C$ and we can
write $P=Q+R$ for $Q\in \Sigma$
and $R\in C$. Then $R$ is rational so $\sigma P-P=\sigma Q-Q$
for any $\sigma \in G$. This implies that $Q$ is also almost
rational. Since the points of $\Sigma$ are cyclotomic,
we have $Q\in \Sigma [3]$ (lemma 3).

The conjunction of the previous two paragraphs
shows that $\JN_{a.r.t.}\subset C\oplus \Sigma [3].$
To check equality, one notes:

-Rational points are almost rational, so points of $C$ are a.r.

-$\Sigma [3]$ consists of almost rational points: This
is because $\Sigma [3]$ is either trivial or isomorphic
to $\mu_3$. It's easy to check that points of $\mu_3$
are almost rational.

-A translate of an a.r. point by a rational point is a.r.

We are done.

{\footnotesize DEPARTMENT OF MATHEMATICS, UNIVERSITY OF CALIFORNIA, BERKELEY,
CA 94720, and DEPARTMENT OF MATHEMATICS, UNIVERSITY OF ARIZONA, TUCSON, AZ 85721}

\end{document}